\documentclass[a4paper,10pt]{amsart}
\usepackage{amsmath}
\usepackage{amssymb}
\usepackage{latexsym}
\usepackage{amsbsy}
\usepackage[all]{xy}
\usepackage{amscd}
\usepackage[dvips]{graphicx}

\newtheorem{Prop}{Proposition}[section]
\newtheorem{Thm}[Prop]{Theorem}
\newtheorem{Lemma}[Prop]{Lemma}
\newtheorem{Cor}[Prop]{Corollary}
\newtheorem{Example}[Prop]{Example}
\newtheorem{Remark}[Prop]{Remark}

\newtheorem{Definition}[Prop]{Definition}

\def\bbn{{\mathbb N}}  \def\bbz{{\mathbb Z}}  \def\bbq{{\mathbb Q}} \def\bb1{{\mathbb 1}}

  \def\bbc{{\mathbb C}}

\def\ra{\rightarrow}

\def\ext{\mbox{Ext}\,} 
\def\dim{\mbox{dim}\,}

\def\uq2{U_q(\hat{sl}_2)}

\def\bb{{\bf b}}

\def\nd{{\noindent}}

\def\mc{{\mathcal{C}}}

\def\md{{\mathcal{D}}}

\def\ue{{\underline{e}}}

\begin{document}

\title[ A $\mathbb{Z}$-basis for the cluster algebra associated to an affine quiver]
{A $\mathbb{Z}$-basis for the cluster algebra associated to an
affine quiver}


\author{Ming Ding, Jie Xiao and Fan Xu}
\address{Department of Mathematical Sciences\\
Tsinghua University\\
Beijing 10084, P.~R.~China} \email{m-ding04@mails.tsinghua.edu.cn
(M.Ding),\ jxiao@math.tsinghua.edu.cn (J.Xiao),\
f-xu04@mails.tsinghua.edu.cn (F.Xu)}

\maketitle


\bigskip

\begin{abstract}
The canonical bases of cluster algebras of finite types and rank 2
are given explicitly in \cite{CK2005} and \cite{SZ} respectively.
In this paper, we will deduce  $\mathbb{Z}$-bases for cluster
algebras for affine types $\widetilde{A}_{n,n},\widetilde{D}$ and
$\widetilde{E}$. Moreover, we give an inductive formula for
computing the multiplication between two generalized cluster
variables associated to objects in a tube.

\end{abstract}

\section{Introduction}
Cluster algebras were introduced by S. Fomin and A. Zelevinsky
\cite{FZ} in order to develop a combinatorial approach to study
problems of total positivity. By the work of \cite{MRZ}, the link
between acyclic cluster algebras and representation theory of
quivers found its general framework in \cite{BMRRT} where the
authors introduced the cluster category. Let $Q$ be an acyclic
quiver with vertex set $Q_0=\{1,2,\cdots, n\}$. Let $A=\bbc Q$ be
the path algebra of $Q$ and we denote by $P_i$ the indecomposable
projective $\bbc Q$-module with the simple top $S_i$ corresponding
to $i\in Q_0$ and $I_i$ the indecomposable injective $\bbc
Q$-module with the simple socle $S_i$. Let $\md^b(Q)$ be the
bounded derived category of $\mathrm{mod} \bbc Q$ with the shift
functor $T$ and the AR-translation $\tau$. The cluster category
associated to $Q$ is the orbit category
$\mathcal{C}(Q):=\md^{b}(Q)/F$ with $F=T\circ\tau^{-1}$. Let
$\mathbb{Q}(x_1,\cdots,x_n)$ be a transcendental extension of
$\mathbb{Q}.$ The Caldero-Chapton map of an acyclic quiver $Q$ is
the map
$$X_?^Q: \mathrm{obj}(\mc(Q))\ra\bbq(x_1,\cdots,x_n)$$ defined in \cite{CC} by
the following rules:
            \begin{enumerate}
                \item if $M$ is an indecomposable $\bbc Q$-module, then
                    $$
                        X_M^Q = \sum_{\textbf e} \chi(\mathrm{Gr}_{\ue}(M)) \prod_{i \in Q_0} x_i^{-\left<\ue, s_i\right>-\left <s_i, \underline{\mathrm{dim}}M - \ue\right
                        >};
                    $$
                \item if $M=TP_i$ is the shift of the projective module associated to $i \in Q_0$, then $$X_M^Q=x_i;$$
                \item for any two objects $M,N$ of $\mathcal C_Q$,
                we have
                    $$X_{M \oplus N}^Q=X_M^QX_N^Q.$$
            \end{enumerate}
Here, we denote by $\left <-,-\right >$ the Euler form on $\bbc
Q$-mod and $Gr_{\ue}(M)$ is the $\ue$-Grassmannian of $M,$ i.e.
the variety of submodules of $M$ with dimension vector $\ue.$ For
any object $M\in \mc(Q)$, $X_M^Q$ will be called the generalized
cluster variable for $M$.

We note that the indecomposable $\bbc Q$-modules and $TP_i$ for
$i\in Q_0$ exhaust the indecomposable objects of the cluster
category $\mc(Q)$:
$$\mathrm{ind}-\mc(Q)=\mathrm{ind}-\bbc Q\sqcup\{TP_i:i\in
Q_0\}.$$ Each object $M$ in $\mc(Q)$ can be uniquely decomposed in
the following way:
$$M=M_0\oplus TP_M$$
where $M_0$ is a module and $P_M$ is a projective  module. The
module $M_0$ can be recovered  using the following homological
functor:
$$H^{0}=\mathrm{Hom}_{\mc(Q)}(\bbc Q,-):\mc(Q)\longrightarrow \bbc Q-\mathrm{mod}$$
Hence we can rewrite $M=M_0\oplus TP_M$ as:
$$M=H^{0}(M)\oplus TP_M.$$

Let $R=(r_{ij})$ be a matrix of size $|Q_0|\times |Q_0|$
satisfying
$$r_{ij}=\dim_{\bbc}\mathrm{Ext}^1(S_i,S_j)$$ for any $i, j\in Q_0$. The Caldero-Chapton map can be reformulated by the following rules \cite{Xu}\cite{Hubery2005}:

\begin{enumerate}
\item
$$
X_{\tau P}=X_{TP}=x^{\underline{\mathrm{dim}}{P/rad P}},
X_{\tau^{-1}I}=X_{T^{-1}I}=x^{\underline{\mathrm{dim}} socI}$$ for
any projective $\bbc Q$-module $P$ and any injective $\bbc
Q$-module $I$; \item
$$ X_{M}=\sum_{\ue}\chi(\mathrm{Gr}_{\ue}(M))x^{\ue R+(\underline{\mathrm{dim}}M-\ue)R^{tr}-
\underline{\mathrm{dim}}M }
$$
where $M$ is a $\bbc Q$-module, $R^{tr}$ is the transpose of the
matrix $R$ and $x^{v}=x^{v_1}_1\cdots x^{v_n}_n$ for
$v=(v_1,\cdots,v_n)\in \mathbb{Z}^{n}$.
\end{enumerate}

Let $\mathcal{AH}(Q)$ be the subalgebra of $\bbq(x_1,\cdots, x_n)$
generated by
$$\{X_M, X_{\tau P}\mid M, P\in\mathrm{ mod}-\bbc Q, P
\mbox{ is projective module}\}.$$

Let  $\mathcal{EH}(Q)$ be the subalgebra of $\mathcal{AH}(Q)$
generated by
 $$\{X_{M}\mid M\in \mathrm{ind}-\mc(Q), \ext^1_{\mc{(Q)}}(M,M)=0\
 \}.$$

In \cite{CK2006}, the authors showed that the Caldero-Chapton map
 induces a one to one correspondence between indecomposable
 objects  in $\mc(Q)$ without self-extension and the cluster
 variables of the cluster algebra $\mathcal{A}(Q).$ Hence, one can
 view $\mathcal{EH}(Q)$ as the cluster algebra for a quiver $Q$.
When $Q$ is a simply laced Dynkin quiver, $\mathcal{AH}(Q)$ is
coincides with $\mathcal{EH}(Q)$. In \cite{CK2005}, the authors
showed a $\bbz$-basis of  $\mathcal{AH}(Q).$ When $Q$ is the
Kronecker quiver, $\mathcal{AH}(Q)$ is also coincides with
$\mathcal{EH}(Q)$. In \cite{CZ}, the authors gave a $\bbz$-basis
of $\mathcal{AH}(Q)$ called the semicanonical basis. When $Q$ is
the quiver of type $\widetilde{D}_4$, $\mathcal{AH}(Q)$ is still
equal to $\mathcal{EH}(Q)$ and a $\bbz$-basis is given in
\cite{DX}(also see Section \ref{d4}).

Now let $Q$ be an affine quiver of type $\widetilde{A}_n,
\widetilde{D}$ or $\widetilde{E}$. There are many references about
the theory of representations of affine quivers, for example, see
\cite{CB} and \cite{DR76}. The main goal of this paper is to give
a $\mathbb{Z}$-basis of $\mathcal{AH}(Q)$ for an affine quiver
$Q$. We will prove the following theorem.
\begin{Thm}\label{30}
Let $Q$ be an affine quiver. A $\mathbb{Z}$-basis of
$\mathcal{AH}(Q)$
 is the following
set:
$$\{X_{L},X_{T\oplus R}|\mathrm{\underline{dim}}(T_1\oplus
R_1)\neq\mathrm{\underline{dim}}(T_2\oplus
R_2),Ext^{1}_{\mc(Q)}(T,R)=0,Ext^{1}_{\mc(Q)}(L,L)=0\}$$ where $L$
is any non-regular exceptional object, $R$ is $0$ or any regular
exceptional module and $T$ is $0$ or any indecomposable regular
module with self-extension  and in addition if $Q\neq
\widetilde{A}_{1,1}$ and $\mathrm{\underline{dim}}(T\oplus
R)=m\delta$ for some $m\in \bbn$, then $R=0$ and $T$ is an
indecomposable module of dimension vector $m\delta$ in a
non-homogeneous tube.
\end{Thm}
Moreover, we give an inductive formula for computing the
multiplication between two generalized cluster variables
associated to objects in a tube in Section 6. Our tools are the
following cluster multiplication theorems.

\begin{Thm}\cite{XiaoXu}\label{XX}
(1) For any $A$-modules $V_{\xi'}$, $V_{\eta'}$we have
$$\hspace{0cm}d^1(\xi', \eta')X_{V_{\xi'}} X_{V_{\eta'}} =\int_{\lambda\neq
\xi'\oplus \eta'} \chi(\mathbb{P}\mathrm{Ext}_{A}^1(V_{\xi'},
V_{\eta'})_{V_{\lambda}})
X_{V_{\lambda}}$$$$+\int_{\gamma,\beta}\chi(\mathbb{P}\mathrm{Hom}_{A}(V_{\eta'},\tau
V_{\xi'})_{V_{\beta}[1]\oplus
V_{\gamma}})X_{V_{\beta}}X_{\tau^{-1} V_{\gamma}}
$$
where $V_{\xi'}$ has no projective direct summand and $d^1(\xi',
\eta')=\mathrm{\dim}_{\bbc}\ext_{A}(V_{\xi'},V_{\eta'}).$

(2) For any $A$-module $V_{\xi'}$ and  $P\in \rho$ is projective
Then
 $$d(\rho, \xi')X_{V_{\xi'}}x^{\underline{\mathrm{dim}}P/radP}=\int_{\delta,\iota'}\chi(\mathbb{P}\mathrm{Hom}_{A}(V_{\xi'},I)_{V_{\delta}[1]\oplus
I'})X_{V_{\delta}}x^{\underline{\mathrm{dim}}\mathrm{soc}I'}$$$$+\int_{\gamma,\rho'}\chi(\mathbb{P}\mathrm{Hom}_{A}(P,V_{\xi'})_{P'[1]\oplus
V_{\gamma}})X_{V_{\gamma}}x^{\underline{\mathrm{dim}}P'/radP'}
 $$
where $I=D\mathrm{Hom}_{A}(P,A),$ $I'\in \iota'$ injective, $P'\in
\rho'$ projective, and $d(\rho, \xi')=\mathrm{dim}_{\bbc}
\mathrm{Hom}_{A}(V_{\rho},V_{\xi'})$.
\end{Thm}

\begin{Thm}\cite{CC}\label{AR formula}
Let $Q$ be an acyclic quiver and  $M$ any indecomposable
non-projective $A$-module, then
$$
X_{M}X_{\tau M}=1+X_{E}
$$
where $E$ is the middle term of the Auslander-Reiten sequence
ending in $M$.
\end{Thm}

\begin{Thm}\cite{CK2006}\label{CK}
Let $Q$ be an acyclic quiver and  $M$,$N$ be any two objects in
$\mathcal{C}(Q)$ such that $dim\ Ext^{1}_{\mathcal{C}(Q)}(M,N)=1$,
then
$$
X_{M}X_{N}=X_{B}+X_{B'}
$$
where $B$ and $B'$ are the unique objects such that there exists
non-split triangles
$$M\longrightarrow B\longrightarrow N\longrightarrow M[1]\ and\ N\longrightarrow B'\longrightarrow M\longrightarrow N[1]$$
\end{Thm}

In \cite{Dupont1}, the author construct  a $\bbz$-basis for a
cluster algebra of type $\widetilde{A}$ referred as the
semicanonical basis. It is interesting to compare it to the
$\mathbb{Z}$-bases in this paper.
\section{Numerators of Laurent expansions in generalized cluster variables}
 In the following, we
will suppose that $Q$ is one of
$\widetilde{A}_{n,n},\widetilde{D}$ and $\widetilde{E}$ with an
orientation where every vertex is a sink or a source. For any
object $M=M_0\oplus (\oplus_{i}s_{i}TP_{i})\in \mathcal{C}(Q)$
where $M_0$ is a module, we extend the dimension vector to the
objects in the cluster category by setting
$$\mathrm{\underline{dim}}(M)=\mathrm{\underline{dim}}(M_0)-(s_{1},\cdots,s_{n}).$$

Let $E_i[n]$ be the indecomposable regular module with quasi-socle
$E_i$ and quasi-length $n$ and $X_0=1$. For any $\bbc Q$-module
$M$, we denoted by $\mathrm{dim}_{\bbc} M(i)$ the $i$-th component
of $\underline{\mathrm{dim}} M$.
\begin{Definition}\label{p}
For $M, N\in \mc(Q)$ with
$\mathrm{\underline{dim}}(M)=(m_{1},\cdots,m_{n})$ and
$\mathrm{\underline{dim}}(N)=(r_{1},\cdots,r_{n})$, we write
$\mathrm{\underline{dim}}(M)\preceq \mathrm{\underline{dim}}(N)$
if $m_{i}\leq r_{i}\ for\ 1\leq i\leq n$. Moreover, if there
exists some i such that $m_{i}< r_{i}$, then we write
$\mathrm{\underline{dim}}(M)\prec \mathrm{\underline{dim}}(N).$

\end{Definition}
\begin{Remark}
Note that for the cluster multiplication formula in the Theorem
\ref{XX}(1), we have the following  exact sequences:
$$0\longrightarrow V_{\eta'}\longrightarrow V_{\lambda} \longrightarrow
V_{\xi'}\longrightarrow  0$$ and
$$0\longrightarrow V_{\beta}\longrightarrow V_{\eta'}\longrightarrow \tau
V_{\xi'}\longrightarrow V_{\gamma}\longrightarrow 0.$$ If
$\mathrm{\dim}_{\bbc}\ext_{A}(V_{\xi'},V_{\eta'})\neq
0\Longrightarrow
\mathrm{\underline{dim}}V_{\beta}+\mathrm{\underline{dim}}\tau^{-1}
V_{\gamma}\prec
\mathrm{\underline{dim}}V_{\eta'}+\mathrm{\underline{dim}}V_{\xi'}=\mathrm{\underline{dim}}V_{\lambda}$
because $\tau V_{\xi'}$ has no injective summand.
\end{Remark}
According to the definition of the Caldero-Chapton map, we
consider the Laurent expansions in generalized cluster variables
$X_{M}=\frac{P(x)}{\Pi_{1\leq i\leq n}x^{m_{i}}_{i}}$ for
$M\in\mathcal{C}(Q)$ such that the integral polynomial $P(x)$ in
the variables $x_i$  is not divisible by any $x_{i}.$ We define
the denominator vector of $X_{M}$ as $(m_{1},\cdots,m_{n})$
\cite{Dupont1}. The following theorem is called as the denominator
theorem.
\begin{Thm}\cite{CK2005}\label{d}
Let Q be an acyclic quiver. Then for any object M in $\mc(Q)$, the
denominator vector of $X_{M}$ is $\mathrm{\underline{dim}}(M).$
\end{Thm}
According to the denominator Theorem \ref{d}, we can prove the
following propositions.
\begin{Prop}\label{1}
If $M$ is $P_i$ or $I_i$ for $1\leq i\leq n$, then
$X_{M}=\frac{P(x)}{x^{\mathrm{\underline{dim}}(M)}}$ where the
constant term of P(x) is 1.
\end{Prop}
\begin{proof}
1) If $i$ is a sink point, we have the following short exact
sequence:
$$0\longrightarrow P_i\longrightarrow I_i\longrightarrow I'\longrightarrow 0$$
Then by the cluster multiplication theorem in Theorem \ref{XX}, we
have:
$$X_{\tau P_i}X_{P_i}=x^{\underline{dim}soc I'}+1$$
Thus the constant term of numerator in  $X_{P_i}$ as an
irreducible fraction of integral polynomials in the variables
$x_i$ is $1$ because of $X_{\tau P_i}=x_i.$

 If $i$ is a source point, we have the following short exact
sequence:
$$0\longrightarrow P'\longrightarrow P_i\longrightarrow I_i\longrightarrow 0$$
Similarly we have:
$$X_{\tau P_i}X_{P_i}=X_{P'}+1$$
Thus we can finish it by induction on $P'$.

2) For $X_{I_i}$, it is totally similar.
\end{proof}
Note that $X_{\tau {P_i}}=x_i=\frac{1}{x^{-1}_i}$ and
$\mathrm{\underline{dim}}(\tau {P_i})=(0,\cdots,0,-1,0,\cdots,0)$
with i-th component 1 and others 0. Hence we denote the
denominator of $X_{\tau {P_i}}$ by $x^{-1}_i$, and assert the
constant term of numerator in $X_{\tau {P_i}}$ is 1. With these
notations, we have the following Proposition \ref{2}.
\begin{Prop}\label{2}
For any object $M \in \mathcal{C}(Q)$, then
$X_{M}=\frac{P(x)}{x^{\mathrm{\underline{dim}}(M)}}$ where the
constant term of P(x) is 1.
\end{Prop}
\begin{proof}
It is enough to consider the case that $M$ is an indecomposable
module.

1) When $M$ is an indecomposable preprojective module, then by
exchange relation in Thereom \ref{AR formula} we have
$$X_{M}X_{\tau M}=\prod_iX_{B_i}+1.$$
Thus by Proposition \ref{1}, we can prove that
$X_{M}=\frac{P(x)}{x^{\mathrm{\underline{dim}}(M)}}$ where the
constant term of $P(x)$ is 1 by induction with the help of the
directness of AR-quiver in the preprojective component of
mod$\bbc Q$. The discussion is similar for any indecomposable
preinjective module.

2) When M is an indecomposable regular module, we only need to
prove that the proposition holds for any regular simple module
according to the exchange relations.

First, suppose M is in some homogeneous tube with dimension vector
$\delta$. Note that there exists a point $e$ such that
$\mathrm{dim}\delta(e)=1.$ Thus we have $$\mathrm{dim} Ext^1_{\bbc
Q}(M,P(e))=\mathrm{dim}Hom_{\bbc Q}(P(e), M)=1.$$ Then we obtain
the following two non-split exact sequences:
$$0\longrightarrow P(e)\longrightarrow L\longrightarrow M\longrightarrow 0$$
and
$$0\longrightarrow L'\longrightarrow P(e)\longrightarrow  M\longrightarrow L''\longrightarrow 0$$
where L and $L'$ are preprojective modules and $L''$  is a
preinjective module. Using Theorem \ref{XX} or Theorem \ref{CK},
we have
$$X_{M}X_{P(e)}=X_{L}+X_{L'}X_{\tau^{-1} L''}$$
where $\mathrm{\underline{dim}}(L'\oplus \tau^{-1} L'')\prec
\mathrm{\underline{dim}}(P(e)\oplus M)$.

 We have already known
that the constant term of the numerator in $X_{P(e)},X_{L}$ as an
irreducible fraction of integral polynomials in the variables
$x_i$ is 1 by 1), then the constant term of the numerator in
$X_{M}$ must be 1.

Now we consider these non-homogeneous tubes. Note that by the AR
formula in Theorem \ref{AR formula}, we only need to prove the
constant term of the numerator in $X_{E_{i}}$ is 1 for $1\leq
i\leq r$. Suppose $M$ is a regular simple module such that
$\mathrm{dim}M(e)=\mathrm{dim}\delta(e)=1.$ We denote $M$ by
$E_{1}$, thus $\mathrm{dim}E_{1}(e)=1$, and $\mathrm{dim}
E_{i}(e)=0$ for $2\leq i\leq r.$ Therefore $\mathrm{dim}
Ext^1(E_{2},P(e))=1$, then we have the following non-split exact
sequences combining the relation $\tau E_{2}=E_{1}$
$$0\longrightarrow P(e)\longrightarrow L\longrightarrow E_{2}\longrightarrow 0$$
and
$$0\longrightarrow L'\longrightarrow P(e)\longrightarrow E_{1}\longrightarrow L''\longrightarrow 0$$
where $L$ and $L'$ are preprojective modules and $L''$  is a
preinjective module. Then  we have
$$X_{E_{2}}X_{P(e)}=X_{L}+X_{L'}X_{\tau^{-1} L''}$$
where $\mathrm{\underline{dim}}(L'\oplus \tau^{-1} L'')\prec
\mathrm{\underline{dim}}(P(e)\oplus E_{2})$.

 We have already known that the
constant term of the numerator in $X_{P(e)},X_{L}$ as an
irreducible fraction of integral polynomials in the variables
$x_i$ is 1 by 1), then the constant term of the numerator in
$X_{E_{2}}$ must be 1.

Note that $\mathrm{dim}E_{1}[2](e)=\mathrm{dim}
E_{1}(e)+\mathrm{dim}E_{2}(e)=1$, by similar discussions, we can
obtain the constant term of the numerator in $X_{E_{1}[2]}$ must
be 1. Thus by $X_{E_{1}}X_{E_{2}}=X_{E_{1}[2]}+1$, we obtain that
the constant term of the numerator in $X_{E_{1}}$ must be 1. Using
the same method, we can prove the constant term of the numerator
in $X_{E_{i}}$ must be 1 for $3\leq i\leq r.$
\end{proof}

\section{Generalized
cluster variables on tubes} Let $Q$ be an affine quiver with the
minimal imaginary root $\delta=(\delta_i)_{i\in Q_0}$. Then tubes
of indecomposable regular $\bbc Q$-modules are indexed by the
projective line $\mathbb{P}^1.$ Let $\lambda$ be the index of a
homogeneous tube and $M(\lambda)$ be the regular simple $\bbc
Q$-module with dimension vector $\delta$ in this homogeneous tube.
Let $M[i]$ be the regular module with quasi-socle $M$ and
quasi-regular length $i$ for any $i\in \bbn$. Let $X_M$ be the
generalized cluster variable associated to $M$ by the
reformulation of the Caldero-Chapton map. Then we have

\begin{Prop}\label{independent}
Let $\lambda$ and $\mu$ be in $\mathbb{P}^1$ such that
$M(\lambda)$ and $M(\mu)$ are two regular simple modules of
dimension vector $\delta$. Then $X_{M(\lambda)}=X_{M(\mu)}.$
\end{Prop}
\begin{proof}
 Choose a vertex $p\in Q_0$ such that
$\delta_p=1.$  We assume that $p$ is a sink. Let $Q$ have the
underlying graph not of type $\widetilde{A}_n.$ Then there is
unique edge $\alpha\in Q_1$ with the head $p$ and tail $p'.$ It is
easy to check $\delta_{p'}=2.$ Let $P_p$ be the indecomposable
projective module corresponding $p$ and $I(\delta-p)$ be the
indecomposable preinjective module of dimension vector
$\delta-\mathrm{\underline{dim}}S_p$. Then
$\mathrm{dim}_{\bbc}Ext^1_{\bbc Q}(I(\delta-p), P_p)=2$. Given any
$\epsilon\in Ext^1_{\bbc Q}(I(\delta-p), P_p),$ we have a short
exact sequence whose equivalence class is $\epsilon$ as follows:
$$
\xymatrix{\varepsilon:\quad 0\ar[r]& P_p\ar[rr]^{\left(%
\begin{array}{c}
  1 \\
  0 \\
\end{array}%
\right)}&&M_{\epsilon}\ar[rr]^{\left(%
\begin{array}{cc}
  0 & 1 \\
\end{array}%
\right)}&&I(\delta-p)\ar[r]&0}
$$
where $(M_{\epsilon})_{i}=(P_p)_i\oplus I(\delta-p)_i$ for any
$i\in Q_0$, $(M_{\epsilon})_{\beta}=I(\delta-p)_{\beta}$ for
$\beta\neq \alpha$ and $(M_{\epsilon})_{\alpha}$ is
$$
(M_{\epsilon})_{\alpha}=\left(%
\begin{array}{cc}
  0 & m(\epsilon, \alpha) \\
  0 & 0 \\
\end{array}%
\right)
$$ where $m(\epsilon, \alpha)\in Hom_{\bbc}(I(\delta-p)_p', P_p).$ For any $\epsilon,
\epsilon'\in Ext^1_{\bbc Q}(I(\delta-p), P_p)$, $M_{\epsilon}\cong
M_{\epsilon'}$ if and only if $m(\epsilon,
\alpha)=tm(\epsilon',\alpha)$ for some $t\in \bbc.$ The regular
simple $\bbc Q$-modules (denoted by $M(\lambda)$) with dimension
vector $\delta$ satisfy that $M(\lambda)_{\alpha}$ is as follows
$$
\xymatrix{M(\lambda)_{p'}=\bbc^2\ar[rr]^{\left(%
\begin{array}{cc}
  1 & \lambda \\
\end{array}%
\right)}&&M(\lambda)_{p}=\bbc}
$$
where $\lambda\in \bbc^*\setminus\{1\}.$  Let $M(\lambda)$ and
$M(\lambda')$ be any two regular simple $\bbc Q$-modules with
dimension vector $\delta$. Let $P$ is an indecomposable projective
module such that $P\subseteq M(\lambda)$ and
$(\mathrm{\underline{dim}}P)_p=0$. Then $P$ is also a submodule of
$M(\lambda')$ and $(\mathrm{\underline{dim}}P)_{p'}=0$. Let $P$ be
an indecomposable projective module such that $P\subseteq
M(\lambda)$ and $(\mathrm{\underline{dim}}P)_p=1$.  Assume that
$P_{\alpha}$ is  $\xymatrix{\bbc\ar[r]^{a+b\lambda}&\bbc}$, then
there exists $P'\in Gr_{\ue}(M(\lambda'))$  such that $P'\cong P$
and $P'_{\alpha}$ is $\xymatrix{\bbc\ar[r]^{a+b\lambda'}&\bbc}$.
Since $\tau M(\lambda)=M(\lambda)$, we know $\tau^{-i}P$ and
$\tau^{-1}P'$ are the submodules of $M(\lambda)$ and $M(\lambda')$
for any $i\in \bbn$, respectively. Hence, any preprojective
submodule $X$ of $M(\lambda)$ corresponds to a preprojective
submodule $X'$ of $M(\lambda')$ and $X\cong X'$.  Let $Q$ be of
type $\widetilde{A}_n$. Then there are two adjacent edge $\alpha$
and $\beta$. Any regular simple module $M(\lambda)$ satisfies that
$M(\lambda)_{\alpha}$ is as follows:
$$
\xymatrix{\bbc\ar[r]^1&\bbc&\bbc\ar[l]_{\lambda}}.$$ The
discussion is similar as above.  If $p$ is a source, the
discussion is also similar. Therefore, there is a homeomorphism
between $Gr_{\ue}(M(\lambda))$ and $Gr_{\ue}(M(\lambda'))$ for any
dimension vector $\ue.$ By definition,
$X_{M(\lambda)}=X_{M(\mu)}.$
\end{proof}

We note that there is an alternative proof of Proposition
\ref{independent} in \cite[Lemma 3.14]{Dupont1}.

\begin{Prop}\label{def1}
For any $m,n\in \bbn$ and $m\geq n$, we have
$$X_{M[m]}X_{M[n]}=X_{M[m+n]}+X_{M[m+n-2]}+\cdots+X_{M[m-n+2]}+X_{M[m-n]}.$$
\end{Prop}

\begin{proof}
When $n=1$, we know
$\mathrm{dim}_{\bbc}\mathrm{Ext}^{1}(M[m],M)=\mathrm{dim}_{\bbc}
\mathrm{Hom}(M, M[m])=1$. The involving non-split short exact
sequences are
$$0\longrightarrow M\longrightarrow M[m+1]\longrightarrow M[m]\longrightarrow 0$$
and
$$0\longrightarrow M\longrightarrow M[m]\longrightarrow M[m-1]\longrightarrow 0.$$
Thus by the cluster multiplication theorem in Theorem \ref{XX} or
Theorem \ref{CK}  and the fact $\tau M[k]=M[k]$ for any $k\in
\bbn$, we obtain the equation
 $$
X_{M[m]}X_{M}=X_{M[m+1]}+X_{M[m-1]}.$$

Suppose that it is right for  $n\leq k$. When $n=k+1,$  we have
$$X_{M[m]}X_{M[k+1]}=X_{M[m]}(X_{M[k]}X_{M}-X_{M[k-1]})= X_{M[m]}X_{M[k]}X_{M}-X_{M[m]}X_{M[k-1]}$$
$$\hspace{-1.3cm}=\sum_{i=0}^{k}X_{M[m+k-2i]}X_{M}-\sum_{i=0}^{k-1}X_{M[m+k-1-2i]}$$
$$\hspace{1.3cm}=\sum_{i=0}^{k}(X_{M[m+k+1-2i]}+X_{M[m+k-1-2i]})-\sum_{i=0}^{k-1}X_{M[m+k-1-2i]}$$
$$
\hspace{-4.6cm}=\sum_{i=0}^{k+1}X_{M[m+k+1-2i]}.$$
\end{proof}
 By
Proposition \ref{independent} and Proposition \ref{def1}, we can
define $X_{n\delta}:=X_{M[n]}$ for $n\in \bbn.$

Now consider the non-homogeneous tubes $T(k)$ with rank $r_k$
 for $1\leq k\leq m$. The corresponding
regular simple modules are $E_{k,1},\cdots,E_{k,r_k}$ with $\tau
E_{k,i+1}=E_{k,i}$. In fact $m=2$ or $3$ in our conditions. If we
restrict the discussion to one tube, we will omit the index $k$
for convenience. Set $q.soc(E_{k,i}[n])=E_{k,i}$ and
$X_{n\delta_{k,i}}=X_{E_{k,i}[nr_k]}$ for $n\in \mathbb{N}$.

\begin{Prop}\label{17}
Let $E_i$ and $E_j$ be two regular simples in a non-homogeneous
tube with rank r. Then we have
$$X_{E_i[mr]}=X_{E_j[mr]}+X_{E_{i+1}[mr-2]}-X_{E_{j+1}[mr-2]}$$
where $1\leq i< j\leq r\ and\ m\in \mathbb{N}.$
\end{Prop}
\begin{proof}
It is easy to prove that
\begin{eqnarray}
  X_{E_i}X_{E_{i+1}[mr-1]}&=& X_{E_{i}[mr]}+X_{E_{i+2}[mr-2]}, \nonumber\\
  X_{E_{i+1}[mr-1]}X_{E_i} &=& X_{E_{i+1}[mr]}+X_{E_{i+1}[mr-2]}.\nonumber
\end{eqnarray}
Hence, we have
$$X_{E_{i}[mr]}=X_{E_{i+1}[mr]}+X_{E_{i+1}[mr-2]}-X_{E_{i+2}[mr-2]}.$$
Similarly we have
\begin{eqnarray}
  X_{E_{i+1}[mr]}&=& X_{E_{i+2}[mr]}+X_{E_{i+2}[mr-2]}-X_{E_{i+3}[mr-2]},\nonumber\\
  X_{E_{i+2}[mr]}&=& X_{E_{i+3}[mr]}+X_{E_{i+3}[mr-2]}-X_{E_{i+4}[mr-2]}, \nonumber\\
  && \vdots \nonumber\\
  X_{E_{j-1}[mr]}&=& X_{E_{j}[mr]}+X_{E_{j}[mr-2]}-X_{E_{j+1}[mr-2]}.\nonumber
\end{eqnarray}
Thus
$$X_{E_{i}[mr]}+X_{E_{i+1}[mr]}+\cdots +X_{E_{j-1}[mr]}$$
$$=(X_{E_{i+1}[mr]}+X_{E_{i+1}[mr-2]}-X_{E_{i+2}[mr-2]})+(X_{E_{i+2}[mr]}+X_{E_{i+2}[mr-2]}-X_{E_{i+3}[mr-2]})$$
$$+\cdots+(X_{E_{j}[mr]}+X_{E_{j}[mr-2]}-X_{E_{j+1}[mr-2]})$$
$$\hspace{-1.8cm}=X_{E_{i+1}[mr]}+X_{E_{i+2}[mr]}+\cdots +X_{E_{j}[mr]}+X_{E_{i+1}[mr-2]}-X_{E_{j+1}[mr-2]}.$$
Therefore
$$X_{E_i[mr]}=X_{E_j[mr]}+X_{E_{i+1}[mr-2]}-X_{E_{j+1}[mr-2]}.$$
\end{proof}
\begin{Example}
Consider r=3, we have
$$X_{E_1}X_{E_2}X_{E_3}=(X_{E_1[2]}+1)X_{E_3}=X_{E_1[3]}+X_{E_1}+X_{E_3}.$$
Similarly we have
$$X_{E_2}X_{E_3}X_{E_1}=X_{E_2[3]}+X_{E_2}+X_{E_1},$$
and
$$X_{E_3}X_{E_1}X_{E_2}=X_{E_2[3]}+X_{E_3}+X_{E_2}.$$
Therefore
$$X_{E_1[3]}=X_{E_2[3]}+X_{E_2}-X_{E_3}=X_{E_3[3]}+X_{E_2}-X_{E_1}.$$
\end{Example}

\begin{Prop}\label{4}
$X_{n\delta_{i,1}}=X_{n\delta_{j,1}}+\sum_{\mathrm{\underline{dim}}L\prec
n\delta}a_LX_L$ for $i\neq j$ in different non-homogeneous tubes
where $a_L\in \mathbb{Q}$.
\end{Prop}
\begin{proof}
Denote $\delta=(v_1,v_2,\cdots,v_n)$ and
$\mathrm{\underline{dim}}(E_{i,j})=(v_{j1},v_{j2},\cdots,v_{jn})$,
then $\delta=\sum_{1\leq j\leq
r_i}\mathrm{\underline{dim}}(E_{i,j})$. Thus by the cluster
multiplication theorem in Theorem \ref{XX}, and the fact that for
any dimension vector there is at most one exceptional module up to
isomorphism, we have
$$\hspace{-8.0cm}X^{v_1}_{S_1}X^{v_2}_{S_2}\cdots X^{v_n}_{S_n}$$
$$\hspace{-4.2cm}=(X^{v_{11}}_{S_1}X^{v_{12}}_{S_2}\cdots X^{v_{1n}}_{S_n})
\cdots (X^{v_{r_i1}}_{S_1}X^{v_{r_i2}}_{S_2}\cdots
X^{v_{r_in}}_{S_n})$$
$$=(a_1X_{E_{i,1}}+\sum_{\mathrm{\underline{dim}}L'\prec
\mathrm{\underline{dim}}E_{i,1}}a_{L'}X_{L'})\cdots
(a_{r_i}X_{E_{i,r_i}}+\sum_{\mathrm{\underline{dim}}L''\prec
\mathrm{\underline{dim}}E_{i,r_i}}a_{L''}X_{L''})$$
$$\hspace{-6.2cm}=a_1\cdots a_{r_i}X_{\delta_{i,1}}+\sum_{\mathrm{\underline{dim}}M\prec
\delta}a_MX_M.$$ Similarly
$$\hspace{-8.0cm}X^{v_1}_{S_1}X^{v_2}_{S_2}\cdots X^{v_n}_{S_n}$$
$$=(b_1X_{E_{j,1}}+\sum_{\mathrm{\underline{dim}}T'\prec
\mathrm{\underline{dim}}E_{j,1}}b_{T'}X_{T'})\cdots
(b_{r_j}X_{E_{j,r_j}}+\sum_{\mathrm{\underline{dim}}T''\prec
\mathrm{\underline{dim}}E_{j,r_j}}b_{T''}X_{T''})$$
$$\hspace{-6.4cm}=b_1\cdots b_{r_j}X_{\delta_{j,1}}+\sum_{\mathrm{\underline{dim}}N\prec
\delta}b_NX_N.$$ Thus we have
$$a_1\cdots a_nX_{\delta_{i,1}}=b_1\cdots b_nX_{\delta_{j,1}}+\sum_{\mathrm{\underline{dim}}N\prec
\delta}b_NX_N-\sum_{\mathrm{\underline{dim}}M\prec
\delta}a_MX_M.$$ Therefore by Proposition \ref{2} and the
denominator theorem in Theorem \ref{d}, we have:
$$X_{\delta_{i,1}}=X_{\delta_{j,1}}+\sum_{\mathrm{\underline{dim}}N'\prec
\delta}b_{N'}X_{N'}.$$ \\
Now, suppose the proposition holds for
$k\leq n,$ then on the one hand
$$X_{n\delta_{i,1}}X_{\delta_{i,1}}=X_{(n+1)\delta_{i,1}}+\sum_{\mathrm{\underline{dim}}L'\prec
(n+1)\delta}b_{L'}X_{L'}.$$ \\
On the other hand by
induction\begin{eqnarray}
  && X_{n\delta_{i,1}}X_{\delta_{i,1}} \nonumber\\
   &=& (X_{n\delta_{j,1}}+\sum_{\mathrm{\underline{dim}}L\prec
n\delta}a_LX_L)X_{\delta_{i,1}} \nonumber\\
   &=& (X_{n\delta_{j,1}}+\sum_{\mathrm{\underline{dim}}L\prec
n\delta}a_LX_L)(X_{\delta_{j,1}}+\sum_{\mathrm{\underline{dim}}N'\prec
\delta}b_{N'}X_{N'})\nonumber \\
   &=& X_{(n+1)\delta_{j,1}}+\sum_{\mathrm{\underline{dim}}L''\prec
(n+1)\delta}b_{L''}X_{L''}.\nonumber
\end{eqnarray}\\
Therefore
$$X_{(n+1)\delta_{i,1}}=X_{(n+1)\delta_{j,1}}+\sum_{\mathrm{\underline{dim}}L''\prec
(n+1)\delta}b_{L''}X_{L''}-\sum_{\mathrm{\underline{dim}}L'\prec
(n+1)\delta}b_{L'}X_{L'}.$$ \\
Thus the proof is finished.
\end{proof}
\begin{Prop}\label{5}
$X_{n\delta}=X_{n\delta_{i,1}}+\sum_{\mathrm{\underline{dim}}L\prec
n\delta}a_LX_L$, where $a_L\in \mathbb{Q}$.
\end{Prop}
\begin{proof}
Suppose $Q$ be of type $\widetilde{D}$ or $\widetilde{E}$ and
$\delta=(v_1,v_2,\cdots,v_n)$, using the same method in
Proposition \ref{4}, we have
$$X^{v_1}_{S_1}X^{v_2}_{S_2}\cdots X^{v_n}_{S_n}$$
$$=(a_1X_{E_{i,1}}+\sum_{\mathrm{\underline{dim}}L'\prec
\mathrm{\underline{dim}}E_{i,1}}a_{L'}X_{L'})\cdots
(a_{r_i}X_{E_{i,r_i}}+\sum_{\mathrm{\underline{dim}}L''\prec
\mathrm{\underline{dim}}E_{i,r_i}}a_{L''}X_{L''})$$
$$\hspace{-6.2cm}=a_1\cdots a_{r_i}X_{\delta_{i,1}}+\sum_{\mathrm{\underline{dim}}M\prec
\delta}a_MX_M.$$
We note that there exists a submodule $L_1$
satisfying the following short exact sequence:
$$0\longrightarrow L_1\longrightarrow \delta\longrightarrow L_2\longrightarrow 0$$
where $L_1,L_2$ are preprojective, preinjective modules
respectively. Thus by Proposition \ref{17} and Proposition
\ref{4}, we have
$$X^{v_1}_{S_1}X^{v_2}_{S_2}\cdots X^{v_n}_{S_n}=(b_1X_{L_{1}}+\sum_{\mathrm{\underline{dim}}L'\prec
\mathrm{\underline{dim}}L_{1}}b_{L'}X_{L'})(b_2X_{L_{2}}+\sum_{\mathrm{\underline{dim}}L''\prec
\mathrm{\underline{dim}}L_{2}}b_{L''}X_{L''})$$
$$=b_1b_2X_{\delta}+\sum^{m}_{k=1}b_{k,1}X_{\delta_{k,1}}+\sum_{\mathrm{\underline{dim}}N\prec
\delta}b_{N}X_{N}$$
$$\hspace{-0.7cm}=b_1b_2X_{\delta}+bX_{\delta_{i,1}}+\sum_{\mathrm{\underline{dim}}N'\prec
\delta}b_{N'}X_{N'}.$$ Thus
$$a_1\cdots a_nX_{\delta_{i,1}}+\sum_{\mathrm{\underline{dim}}M\prec
\delta}a_MX_M=b_1b_2X_{\delta}+bX_{\delta_{i,1}}+\sum_{\mathrm{\underline{dim}}N'\prec
\delta}b_{N'}X_{N'}.$$ Therefore by Proposition \ref{2} and the
denominator theorem in \cite{Hubery2005}, we have
$$X_{\delta}=X_{\delta_{i,1}}+\sum_{\mathrm{\underline{dim}}M'\prec
\delta}a_{M'}X_{M'}.$$ Then we can finish the proof by induction
as in the proof of Proposition \ref{4}.

Now we assume $Q$ is of the form $\widetilde{A}$. We give an
alternative proof of the difference property in \cite{Dupont1}.
Let $Q$ be a quiver as follows \cite{DR76}:
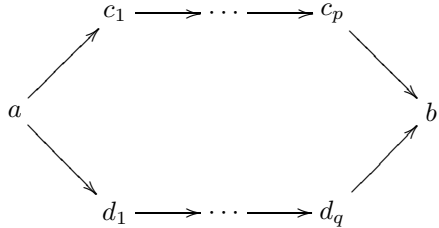
$$
\xymatrix{&c_1\ar[r]&\cdots\ar[r]&c_p\ar[rd]&\\
a\ar[ur]\ar[dr]&&&&b\\
&d_1\ar[r]&\cdots\ar[r]&d_q\ar[ur]&}
$$
Let $\lambda\in \bbc^{*}$ and $M(\lambda)$ be the regular simple
module as follows
$$
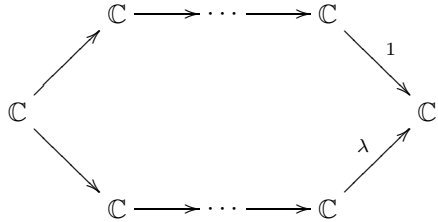
\xymatrix{&\bbc\ar[r]&\cdots\ar[r]&\bbc\ar[rd]^{1}&\\
\bbc\ar[ur]\ar[dr]&&&&\bbc\\
&\bbc\ar[r]&\cdots\ar[r]&\bbc\ar[ur]^{\lambda}&}
$$
Its proper submodules $M_0(\lambda)$ are of the forms as follows:
$$
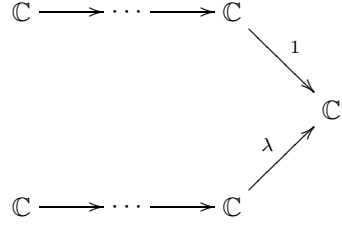
\xymatrix{&\bbc\ar[r]&\cdots\ar[r]&\bbc\ar[rd]^{1}&\\
&&&&\bbc\\
&\bbc\ar[r]&\cdots\ar[r]&\bbc\ar[ur]^{\lambda}&}
$$
Let $M(0)$ be the regular module as follows
$$
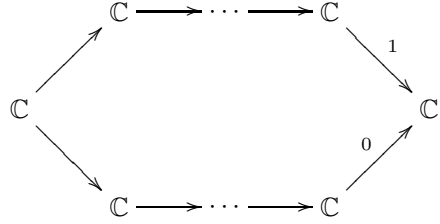
\xymatrix{&\bbc\ar[r]&\cdots\ar[r]&\bbc\ar[rd]^{1}&\\
\bbc\ar[ur]\ar[dr]&&&&\bbc\\
&\bbc\ar[r]&\cdots\ar[r]&\bbc\ar[ur]^{0}&}
$$
Its proper submodules $M_0(0)$ are of the following two forms:
$$
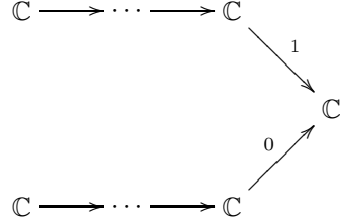
\xymatrix{&\bbc\ar[r]&\cdots\ar[r]&\bbc\ar[rd]^{1}&\\
&&&&\bbc\\
&\bbc\ar[r]&\cdots\ar[r]&\bbc\ar[ur]^{0}&}
$$
and
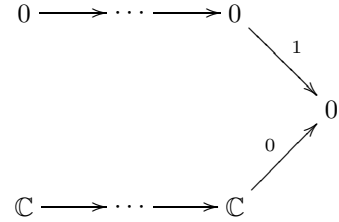
$$
\xymatrix{&0\ar[r]&\cdots\ar[r]&0\ar[rd]^{1}&\\
&&&&0\\
&\bbc\ar[r]&\cdots\ar[r]&\bbc\ar[ur]^{0}&}
$$
The proper submodules with the second form lie in the
non-homogeneous tube indexed by $0$ with quasi socle $S_{d_{q}}.$
Hence, there is a bijection between the submodules of $M(\lambda)$
and the submodules of $M(0)$ of the first form. It is easy to
conclude the difference property in \cite{Dupont1}:
$$X_{M(0)}=X_{M(\lambda)}+X_{q.rad M(0)/S_{d_{q}}}.$$

\end{proof}

\begin{Prop}\label{6}
If $\mathrm{\underline{dim}}(T_1\oplus
R_1)=\mathrm{\underline{dim}}(T_2\oplus R_2)$ where $R_i$ are 0 or
any regular exceptional modules, $T_i$ are 0 or any indecomposable
regular modules with self-extension in non-homogeneous tubes and
there are no extension between $R_i$ and $T_i$, then
$$X_{T_1\oplus R_1}=X_{T_2\oplus R_2}+\sum_{\mathrm{\underline{dim}}R\prec
\mathrm{\underline{dim}}(T_2\oplus R_2)}a_{R}X_{R}$$ where
$a_{R}\in \mathbb{Q}$.
\end{Prop}
\begin{proof}
Suppose $\mathrm{\underline{dim}}(T_1\oplus
R_1)=(d_1,d_2,\cdots,d_n)$, using the same method in Proposition
\ref{4}, we have
$$X^{d_1}_{S_1}X^{d_2}_{S_2}\cdots X^{d_n}_{S_n}=(a_1X_{E_{i,1}}+\sum_{\mathrm{\underline{dim}}L'\prec
\mathrm{\underline{dim}}E_{i,1}}a_{L'}X_{L'})\cdots
(a_sX_{E_{i,s}}+\sum_{\mathrm{\underline{dim}}L''\prec
\mathrm{\underline{dim}}E_{i,s}}a_{L''}X_{L''})$$
$$\times(a_{R_1}X_{R_1}+\sum_{\mathrm{\underline{dim}}L'''\prec
\mathrm{\underline{dim}}R_1}a_{L'''}X_{L'''})$$
$$=aX_{T_1\oplus R_1}+\sum_{\mathrm{\underline{dim}}L\prec
\mathrm{\underline{dim}}(T_1\oplus R_1)}a_{L}X_{L}.$$

$$X^{d_1}_{S_1}X^{d_2}_{S_2}\cdots X^{d_n}_{S_n}=(b_1X_{E_{j,1}}+\sum_{\mathrm{\underline{dim}}M'\prec
\mathrm{\underline{dim}}E_{j,1}}b_{M'}X_{M'})\cdots
(b_tX_{E_{j,t}}+\sum_{\mathrm{\underline{dim}}M''\prec
\mathrm{\underline{dim}}E_{i,t}}b_{M''}X_{M''})$$
$$\times(b_{R_2}X_{R_2}+\sum_{\mathrm{\underline{dim}}M'''\prec
\mathrm{\underline{dim}}R_2}b_{M'''}X_{M'''})$$
$$=bX_{T_2\oplus R_2}+\sum_{\mathrm{\underline{dim}}L\prec
\mathrm{\underline{dim}}(T_2\oplus R_2)}b_{M}X_{M}.$$ Thus

$$aX_{T_1\oplus R_1}+\sum_{\mathrm{\underline{dim}}L\prec
\mathrm{\underline{dim}}(T_1\oplus R_1)}a_{L}X_{L}=bX_{T_2\oplus
R_2}+\sum_{\mathrm{\underline{dim}}L\prec
\mathrm{\underline{dim}}(T_2\oplus R_2)}b_{M}X_{M}.$$
 Therefore by Proposition \ref{2} and the
denominator theorem in Theorem \ref{d}, we have
$$X_{T_1\oplus R_1}=X_{T_2\oplus R_2}+\sum_{\mathrm{\underline{dim}}R\prec
\mathrm{\underline{dim}}(T_2\oplus R_2)}a_{R}X_{R}.$$
\end{proof}

\section{A $\mathbb{Z}$-basis for the cluster algebra of the alternating
quiver of $\widetilde{A}_{n,n},\widetilde{D}\ or\ \widetilde{E}$}

Recall that for an acyclic quiver, the matrix B associated to Q is
the anti-symmetric matrix given by
$$b_{ij}=|{i\longrightarrow j\in Q_{1}}|-|{j\longrightarrow i\in Q_{1}}|$$
where $1\leq i,j\leq n$.

\begin{Definition}\cite{CK2005}\cite{Dupont}
Let Q be an acyclic quiver with associated matrix B. Q is called
graded if there exists a linear form $\epsilon$ on
$\mathbb{Z}^{n}$ such that $\epsilon(B\alpha_i)<0$ for any $1\leq
i\leq n$ where $\alpha_i$ denotes the i-th vector of the canonical
basis of $\mathbb{Z}^{n}$.
\end{Definition}

\begin{Thm}\cite{CK2005}\label{40}
Let Q be a graded quiver and $\{M_1,\cdots,M_r\}$ a family objects
in $\mc(Q)$ such that $\mathrm{\underline{dim}}(M_i)\neq
\mathrm{\underline{dim}}(M_j)$ for $i\neq j$, then
$X_{M_1},\cdots,X_{M_r}$ are linearly independent over
$\mathbb{Q}$.
\end{Thm}
In the section, we still suppose that $Q$ is one of
$\widetilde{A}_{n,n},\widetilde{D}$ and $\widetilde{E}$ with an
orientation where every vertex is a sink or a source. Note that
the quiver $Q$ we consider is graded. We will prove the following
theorem.
\begin{Thm}\label{7}
A  $\mathbb{Z}$-basis for $\mathcal{AH}(Q)$ is the following set
denoted by $\mathcal{S}(Q)$:
$$\{X_{L},X_{T\oplus R}|\mathrm{\underline{dim}}(T_1\oplus
R_1)\neq\mathrm{\underline{dim}}(T_2\oplus
R_2),Ext^{1}_{\mc(Q)}(T,R)=0,Ext^{1}_{\mc(Q)}(L,L)=0\}$$ where $L$
is any non-regular exceptional object, $R$ is $0$ or any regular
exceptional module and $T$ is $0$ or any indecomposable regular
module with self-extension.

\end{Thm}
\begin{proof}
Note that for any objects $M,N \in \mathcal{C}(Q),$ the final
results of their multiplication $X_{M}X_{N}$ must be
$\mathbb{Q}-$combinatorics of $$\{X_{L},X_{T\oplus
R}|Ext^{1}_{\mc(Q)}(T,R)=0,Ext^{1}_{\mc(Q)}(L,L)=0\}$$ where $L$
is any non-regular exceptional object, $R$ is $0$ or any regular
exceptional module and $T$ is $0$ or any indecomposable regular
module with self-extension.

Hence from those Propositions in Section 3, we can easily find
that $X_{M}X_{N}$ is a $\mathbb{Q}-$combinatorics of
$\mathcal{S}(Q)$. Then the proof is finished by the following
Proposition \ref{9} and Proposition \ref{11}.
\end{proof}
\begin{Remark}\label{15}
By  Proposition \ref{12} in the following, we can also rewrite the
$\mathbb{Z}$-basis as $\{X_{M(d_1,\cdots,d_n)}:(d_1,\cdots,d_n)\in
\mathbb{Z}^n\}$, where $M(d_1,\cdots,d_n)=L\ or\ T\oplus R$
associated to the element in $\mathcal{S}(Q)$.
\end{Remark}

\begin{Cor}\label{81}
$\mathcal{S}(Q)$ is a $\mathbb{Z}$-basis of the cluster algebra
$\mathcal{EH}(Q)$.
\end{Cor}
\begin{proof}
It is obvious that  $X_{E_i}\in \mathcal{EH}(Q)$. Then by
Proposition 6.2 in \cite{Dupont}, one can show that $X_{E_i[n]}\in
\mathcal{EH}(Q)$ for any $n\in \bbn$. Thus by Proposition \ref{5},
we can prove that $X_{m\delta}\in \mathcal{EH}(Q)$. By definition,
$X_{L}\in \mathcal{EH}(Q)$ for $L$ satisfying
$\mathrm{Ext}_{\mathcal{C}(Q)}^1(L,L)=0$. Thus
$\mathcal{S}(Q)\subset \mathcal{EH}(Q)$. It follows that
$\mathcal{S}(Q)$ is a $\mathbb{Z}$-basis of the cluster algebra
$\mathcal{EH}(Q)$ by Theorem \ref{7}.
\end{proof}
According to Theorem \ref{7} and Corollary \ref{81}, we have
\begin{Cor}
$\mathcal{EH}(Q)=\mathcal{AH}(Q)$.
\end{Cor}

Firstly, by Theorem \ref{40} we  need to prove the dimension
vectors of these  objects associated to the corresponding elements
in $\mathcal{S}(Q)$ are different.

\begin{Prop}\label{9}
Let $M$ be a  regular module associated to some element in
$\mathcal{S}(Q)$ and $L$ be a non-regular exceptional object in
$\mc(Q)$. Then $\underline{\mathrm{dim}}(M)\neq
\underline{\mathrm{dim}}(L)$.

\end{Prop}

\begin{proof}
If $L$ contains some $\tau P_i$ as its direct summand, we know
that
$$\underline{\mathrm{dim}}(\tau P_i)=(0,\cdots,0,-1,0,\cdots,0)$$
where the $i$-th component is $-1$. Suppose $L=\tau P_i\oplus\tau
P_{i_1}\cdots\tau P_{i_r}\oplus N$ where $N$ is an exceptional
module. Because  $L$ is an exceptional object, $ X_{\tau
P_i}X_{N}=X_{\tau P_i\oplus N}$ i.e.
$\mathrm{dim}_{\bbc}\mathrm{Hom}(P_i,N)=0$. Thus we have
$\mathrm{dim}_{\bbc} N(i)=0$ and $\mathrm{dim}_{\bbc}(\tau
P_i\oplus\tau P_{i_1}\cdots\tau P_{i_r}\oplus N)(i)\leq -1$.
However, $\underline{\mathrm{dim}}(M)\geq 0$. Therefore,
$\underline{\mathrm{dim}}(M)\neq \underline{\mathrm{dim}}(L)$.

If $L$ is a module. Suppose $\underline{\mathrm{dim}}(M)=
\underline{\mathrm{dim}}(L)$. Because $L$ is an exceptional
module, we know that $M$ belongs to the orbit of $L$ and then $M$
is a degeneration of $L$. Hence, there exists some $\bbc Q$-module
$U$ such that
$$ 0\longrightarrow
M\longrightarrow L\oplus U\longrightarrow U\longrightarrow 0$$ is
an exact sequence. Choose minimal $U$ so that we cannot separate
the following exact sequence
$$0\longrightarrow  0\longrightarrow U_1\longrightarrow U_1\longrightarrow
0$$ from the above short exact sequence. Thus  $M$ has a non-zero
map to every direct summand of $L$. Therefore $L$  has no
preprojective modules as direct summand because  $M$ is a regular
module.

Dually there exists a $\bbc Q$-module $V$ such that $$
0\longrightarrow V\longrightarrow V\oplus L\longrightarrow
M\longrightarrow 0$$ is an exact sequence. We can choose minimal
$V$ so that one cannot separate the following exact sequence
$$0\longrightarrow V_1\longrightarrow V_1\longrightarrow
0\longrightarrow 0$$ from the above short exact sequence. Thus
every direct summand of  $L$ has a non-zero map to $M$. Therefore
$L$  has no preinjective modules as direct summand because  M is a
regular module.

Therefore $L$ is a regular exceptional module, it is a
contradiction.
\end{proof}

Secondly, we need to prove that $\mathcal{S}(Q)$ is a
$\mathbb{Z}-$basis.

\begin{Prop}\label{11}
$X_{M}X_{N}$  belongs to   $\mathbb{Z}\mathcal{S}(Q)$ for any $M,N
\in \mathcal{C}(Q).$
\end{Prop}

\begin{proof}
According to these above discussions, for any objects $M,N \in
\mathcal{C}(Q),$ we have
$$X_{M}X_{N}=b_LX_{L}+\sum_{\mathrm{\underline{dim}}L'\prec
\mathrm{\underline{dim}}(M\oplus N)}b_{L'}X_{L'}$$ where
$\mathrm{\underline{dim}}(L)=\mathrm{\underline{dim}}(M\oplus N)$,
$X_{L},X_{L'}\in \mathcal{S}(Q)$ and $b_L,b_{L'}\in \mathbb{Q}$.

Therefore by Proposition \ref{2} and the denominator theorem in
Theorem \ref{d}, we have $b_L=1$. Note that there exists a partial
order on these dimension vectors by Definition \ref{p}. Thus in
these remained $L'$, we choose these maximal elements denoted by
$L'_{1},\cdots,L'_{s}$. Then by $b_L=1$ and the coefficients of
Laurent expansions in generalized cluster variables are integers,
we obtain that $a_{L'_{1}},\cdots,a_{L'_{s}}$ are integers. Using
the same method, we have $b_{L'}\in \mathbb{Z}$.
\end{proof}

We denote $X^{-d}_{S_i}=X^{d}_{\tau P_i}$ for $d\in \mathbb{N}$.
Then we have the following result.

\begin{Prop}\label{12}
$$X^{d_1}_{S_1}X^{d_2}_{S_2}\cdots X^{d_n}_{S_n}=X_{M(d_1,\cdots,d_n)}+\sum_{\underline{dim}L\prec
(d_1,\cdots,d_n)}b_{L}X_{L}$$ where
$X_{M(d_1,\cdots,d_n)},X_{L}\in \mathcal{S}(Q)$,
$\mathrm{\underline{dim}}M=(d_1,\cdots,d_n)\in \mathbb{Z}^n$ and
$b_{L}\in \mathbb{Z}$.

\end{Prop}

\begin{proof}
By Proposition \ref{2}, Theorem \ref{d} and Theorem \ref{7}.
\end{proof}

Note that $\{X_{M(d_1,\cdots,d_n)}:(d_1,\cdots,d_n)\in
\mathbb{Z}^n\}$ is a $\mathbb{Z}$-basis of $\mathcal{AH}(Q)$, then
we have the following Corollary \ref{13} by Proposition \ref{12}.

\begin{Cor}\label{13}
$\{X^{d_1}_{S_1}X^{d_2}_{S_2}\cdots
X^{d_n}_{S_n}:(d_1,\cdots,d_n)\in \mathbb{Z}^n\}$ is a
$\mathbb{Z}$-basis of $\mathcal{AH}(Q)$, further, there is an
isomorphism of $\mathbb{Z}-$algebras:
$$\mathcal{AH}(Q)\simeq \mathbb{Z}[X_{S_1},\cdots,X_{S_n},X_{\tau P_1},\cdots,X_{\tau P_n}].$$
\end{Cor}

\begin{Remark}\label{14}
In \cite{Dupont}, G.Dupont has proved the following isomorphism of
$\mathbb{Z}-$algebras:
$$\mathcal{A}(\overrightarrow{A_n})\simeq \mathbb{Z}[X_{S_1},\cdots,X_{S_n},X_{\tau P_1},\cdots,X_{\tau P_n}]$$
where $n>0.$
\end{Remark}
\section{A $\mathbb{Z}$-basis for the cluster algebra of type
$\widetilde{A}_{n,n},\widetilde{D}\ or\ \widetilde{E}$}

In this section, we will deduce a $\mathbb{Z}$-basis for the
cluster algebra of $\widetilde{A}_{n,n},\widetilde{D}\ or\
\widetilde{E}$ which is independent of their orientations.

Define the reflected quiver $\sigma_{i}(Q)$ by reversing all the
arrows ending at $i$. The mutations can be viewed as
generalizations of reflections i.e. if $i$ is a sink or a source
in $Q_0,$ then $\mu_{i}(Q)=\sigma_{i}(Q)$ where $\mu_{i}$ denotes
the mutation in the direction $i$. Thus there is a natural
isomorphism of cluster algebras
$$\Phi: \mathcal{A}(Q)\longrightarrow \mathcal{A}(Q')$$
where $Q'$ is a quiver mutation equivalent to $Q$, and $\Phi$ is
called the canonical cluster algebras isomorphism in
\cite{Dupont1}.

Now suppose $Q$ is an acyclic quiver and $i$ is a sink in $Q_0$.
Let $Q'=\sigma_{i}(Q)$ and
$R^{+}_{i}:\mathcal{C}({Q})\longrightarrow \mathcal{C}({Q'})$ be
the extended BGP-reflection functor defined in \cite{Zhu}. Denote
by $X^{Q}_{?}$ (resp. by $X^{\sigma_iQ}_{?}$) the Caldero-Chapton
map associated to $Q$ (resp. to $\sigma_iQ$).

Then the following hold.
\begin{Lemma}\cite{Zhu}\label{80}
Let $Q$ be an acyclic quiver and  $i$ be a sink in $Q$. Then
$R^{+}_{i}$ induces a triangle equivalence
 $$R^{+}_{i}:\mathcal{C}({Q})\longrightarrow
\mathcal{C}({\sigma_iQ})$$

\end{Lemma}

\begin{Lemma}\cite{Dupont1}\label{70}
Let $Q$ be an affine quiver and  $i$ be a sink in $Q$ and $M$ be
an indecomposable regular $\mathbb{C}Q$-module with
quasi-composition series
$$0=M_0\subset M_1\subset \cdots \subset M_s=M$$
then $R^{+}_{i}M$ is a regular $\mathbb{C}\sigma_iQ$-module with
quasi-composition series
$$0=R^{+}_{i}M_0\subset R^{+}_{i}M_1\subset \cdots \subset R^{+}_{i}M_s=R^{+}_{i}M$$
In particular, $R^{+}_{i}$ send quasi-socles to quasi-socles,
quasi-radicals to quasi-radicals and preserves quasi-lengths.
\end{Lemma}

\begin{Lemma}\cite{Dupont1}\label{50}
Let $Q$ be an affine quiver and  $i$ be a sink in $Q$. Denote by
$\Phi_i:\mathcal{A}(Q)\longrightarrow \mathcal{A}(\sigma_iQ)$ the
canonical cluster algebra isomorphism and by
$R^{+}_{i}:\mathcal{C}({Q})\longrightarrow
\mathcal{C}({\sigma_iQ})$ the extended BGP-reflection functor.
Then
$$\Phi_i(X^{Q}_{M})=X^{\sigma_iQ}_{R^{+}_{i}M}$$
where $M$ is any rigid object in $\mathcal{C}({Q})$ or any regular
module in non-homogeneous tubes of rank $r> 1.$
\end{Lemma}
In the following, we will suppose that $Q$ is one of
$\widetilde{A}_{n,n},\widetilde{D}$ and $\widetilde{E}$ with an
orientation where every vertex is a sink or a source, and $Q'$ is
an acyclic quiver of $\widetilde{A}_{n,n},\widetilde{D}$ or
$\widetilde{E}$. Then there exists an admissible sequence of sinks
$(i_1,\cdots,i_s)$ such that $Q'=\sigma_{i_s}\cdots
\sigma_{i_1}(Q)$. Denote by $\Phi:\mathcal{A}(Q)\longrightarrow
\mathcal{A}(Q')$ the canonical cluster algebra isomorphism and by
$R^{+}=R^{+}_{i_s}\cdots
R^{+}_{i_1}:\mathcal{C}({Q})\longrightarrow \mathcal{C}({Q'})$.
Note that $R^{+}$ is an equivalence of triangulated categories by
Lemma \ref{80}. Then we have the following main theorem.

\begin{Thm}\label{60}
A $\mathbb{Z}$-basis for the cluster algebra of $Q'$ is the
following set denoted by $\mathcal{B}(Q')$:
$$\{X_{L'},X_{T'\oplus R'}|\mathrm{\underline{dim}}(T'_1\oplus
R'_1)\neq\mathrm{\underline{dim}}(T'_2\oplus
R'_2),Ext^{1}_{\mc(Q')}(T',R')=0,Ext^{1}_{\mc(Q')}(L',L')=0\}$$
where $L'$ is any non-regular exceptional object, $R'$ is $0$ or
any regular exceptional module, $T'$ is $0$ or any indecomposable
regular module with self-extension and in addition if $Q'\neq
\widetilde{A}_{1,1}$ and $\mathrm{\underline{dim}}(T'\oplus
R')=m\delta$ for some $m\in \bbn$, then $R'=0$ and $T'$ is an
indecomposable module of dimension vector $m\delta$ in a
non-homogeneous tube.
\end{Thm}
\begin{proof}
If $Q'= \widetilde{A}_{1,1}$, it is obvious that
$$\{X_M,X_{n\delta}: M\in
obj\mathcal{C}(Q),Ext^1(M,M)=0\}$$ is a $\mathbb{Z}$-basis for
cluster algebra of $\widetilde{A}_{1,1}$, which is called the
semicanonical basis in \cite{CZ}.

 If $Q'\neq \widetilde{A}_{1,1}$, then by Theorem \ref{7}, we have
 obtained a $\mathbb{Z}$-basis for cluster
algebra of $Q$:
$$\{X_{L},X_{T\oplus R}|\mathrm{\underline{dim}}(T_1\oplus
R_1)\neq\mathrm{\underline{dim}}(T_2\oplus
R_2),Ext^{1}_{\mc(Q)}(T,R)=0,Ext^{1}_{\mc(Q)}(L,L)=0\}$$ where $L$
is any  non-regular exceptional object, $R$ is $0$ or any regular
exceptional module and $T$ is $0$ or any indecomposable regular
module with self-extension.

If $\mathrm{\underline{dim}}(T\oplus R)=m\delta$ for some $m\in
\bbn$, then we take $R=0$ and $T$  an indecomposable module of
dimension vector $m\delta$ in a non-homogeneous tube. We denote
this $\mathbb{Z}$-basis by $\mathcal{B}(Q)$. Thus
$\Phi(\mathcal{B}(Q))$ is a $\mathbb{Z}$-basis for the cluster
algebra of $Q'$ because $\Phi:\mathcal{A}(Q)\longrightarrow
\mathcal{A}(Q')$ is the canonical cluster algebras isomorphism.
Then by Lemma \ref{80}, Lemma \ref{70} and Lemma \ref{50}, we know
that $\Phi(\mathcal{B}(Q))$ is exactly the basis $\mathcal{B}(Q')$
in Theorem \ref{60}.
\end{proof}
\section{some special cases}
In this section, we will give some examples for special cases to
explain the $\mathbb{Z}$-basis explicitly.
\subsection{$\bbz$-basis for finite type}
For finite type, we know that there are no regular parts in
$\mathcal{B}(Q)$. Thus the $\mathbb{Z}$-bases are exactly
$$\{X_M: M\in obj\mathcal{C}(Q),Ext^1(M,M)=0\}$$ which coincides
with the canonical basis for a cluster algebra of finite type in
\cite{CK2005}.
\subsection{$\bbz$-basis for the Kronecker quiver}
Consider the Kronecker quiver, we note that $R=0$ and $T=n\delta$
for $n\geq 1$. Thus a $\mathbb{Z}$-basis is
$$\{X_M,X_{n\delta}: M\in obj\mathcal{C}(Q),Ext^1(M,M)=0\}$$ which
is called the semicanonical basis in \cite{CZ}. If we modify:
$$z_1:=X_{\delta},\ z_n:=X_{n\delta}-X_{(n-2)\delta}.$$
Then $\{X_M,z_n: M\in obj\mathcal{C}(Q),Ext^1(M,M)=0\}$ is the
canonical basis for cluster algebra of Kronecker quiver in
\cite{SZ}.
\subsection{$\bbz$-basis for $\widetilde{D}_4$}\label{d4}
Let $Q$ be the tame quiver of type $\widetilde{D}_4$ as follows
$$
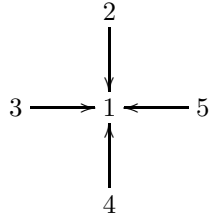
\xymatrix{& 2 \ar[d] &\\
3 \ar[r] & 1  & 5 \ar[l]\\
& 4 \ar[u] &}
$$
We denote the minimal imaginary root by $\delta=(2,1,1,1,1)$. The
regular simple modules of dimension vector $\delta$ are
$$
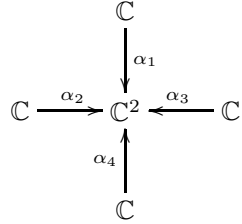
\xymatrix{& \bbc \ar[d]^{\alpha_1} &\\
\bbc \ar[r]^{\alpha_2} & \bbc^2  & \bbc \ar[l]_{\alpha_3}\\
& \bbc \ar[u]^{\alpha_4} &}
$$
with linear maps
$$
\alpha_1=\left(%
\begin{array}{c}
  1 \\
  0 \\
\end{array}%
\right), \alpha_2=\left(%
\begin{array}{c}
  0 \\
  1 \\
\end{array}%
\right),\alpha_3=\left(%
\begin{array}{c}
  1 \\
  1 \\
\end{array}%
\right),\alpha_4=\left(%
\begin{array}{c}
  \lambda\\
  \mu \\
\end{array}%
\right)
$$
where $\lambda/\mu\in \mathbb{P}^1, \lambda/\mu\neq 0,1,\infty.$
Let $M$ be any regular simple $\bbc Q$-module of dimension vector
$\delta$ and $M[i]$ be the regular module with top $M$ and regular
length $i$ for any $i\in \bbn$. Let $X_M$ be the generalized
cluster variable associated to $M$ by the reformulation of the
Caldero-Chapton map. Then we have
\begin{Prop}\cite{DX}
$X_{M}=\frac{1}{x^2_1x_2x_3x_4x_5}+\frac{4}{x_1x_2x_3x_4x_5}+
\frac{x^2_1+4x_1+6}{x_2x_3x_4x_5}+\frac{x_2x_3x_4x_5+2}{x^2_1}+\frac{4}{x_1}.$
\end{Prop}
We define $X_{n\delta}:=X_{M[n]}$ for $n\in \bbn.$ Now, we
consider three non-homogeneous tubes labelled by the subset
$\{0,1,\infty\}$ of $\mathbb{P}^1$. Let
$X_{\delta_1},X_{\delta_2},X_{\delta_3},X_{\delta_4},X_{\delta_5},X_{\delta_6}$
be generalized cluster variables associated to regular modules in
non-homogeneous tubes of dimension vector $\delta$. The regular
simple modules in non-homogeneous tubes are denoted by
$E_1,E_2,E_3,E_4,E_5,E_6$, where
$$\underline{\mathrm{dim}}(E_1)=(1,1,1,0,0),\,\,\underline{\mathrm{dim}}(E_2)=(1,0,0,1,1),\,\,
\underline{\mathrm{dim}}(E_3)=(1,1,0,1,0),$$$$
\underline{\mathrm{dim}}(E_4)=(1,0,1,0,1),\,\,
\underline{\mathrm{dim}}(E_5)=(1,0,1,1,0),\,\,\underline{\mathrm{dim}}(E_6)=(1,1,0,0,1).$$
We note that $\{E_1, E_2\}$, $\{E_3, E_4\}$ and $\{E_5, E_6\}$ are
pairs of the regular simple modules in the bottom of
non-homogeneous tubes labelled by $1,\infty$ and $0$,
respectively.

Let $E_i[n]$ be the indecomposable regular module with top $E_i$
and regular length $n$ for $1\leq i\leq 6$. We set
$X_{n\delta_i}:= X_{E_i[2n]}$ for  $1\leq i\leq 6$.
\begin{Prop}\cite{DX}\label{d4-delta}
For any $n\in \bbn$, we have
$X_{n\delta_1}=X_{n\delta_2}=X_{n\delta_3}=X_{n\delta_4}=X_{n\delta_5}=X_{n\delta_6}$
and $X_{n\delta_1}=X_{n\delta}+X_{(n-1)\delta}$ where $X_0=1$.
\end{Prop}

\begin{Thm}\cite{DX}
A $\mathbb{Z}$-basis for cluster algebra of $\widetilde{D_4}$ is
the following set denoted by $\mathcal{B}(Q)$:
$$\{
X_{L},X_{m\delta},X_{E_i[2k+1]\oplus
R}|\mathrm{\underline{dim}}(E_i[2k+1]\oplus
R_1)\neq\mathrm{\underline{dim}}(E_j[2l+1]\oplus R_2),$$
$$\mathrm{Ext}_{\mathcal{C}(Q)}^1(L,L)=0,\mathrm{Ext}_{\mathcal{C}(Q)}^1(E_i[2k+1],R)=0,m,k,l\geq
0,1\leq i,j\leq 6\}$$ where $L$ is any non-regular exceptional
object, $R$ is $0$ or any regular exceptional module.
\end{Thm}

\section{the inductive multiplication formula for a tube}
Now we fix a tube with rank $r$ and these regular simple modules
are $E_1,\cdots,E_r$ with $\tau E_{i+1}=E_{i}$ where
$E_i=E_{i+mr}$ for $1\leq i\leq r$ and $m\in \mathbb{N}$. Let
$X_{E_i}$ be the corresponding generalized cluster variable for
$i=1,\cdots,r.$ With these notations, we have the following
inductive cluster multiplication formula.
\begin{Thm}\label{16}
Let $i,j, k,l,m$ and $r$ be in $\bbz$ such that $1\leq k\leq
mr+l,0\leq l\leq r-1,1\leq i,j\leq r,m\geq 0$.

\nd (1)When $j\leq i$, then

1)for $k+i\geq r+j$, we have
$X_{E_i[k]}X_{E_j[mr+l]}=X_{E_i[(m+1)r+l+j-i]}X_{E_j[k+i-r-j]}+X_{E_i[r+j-i-1]}X_{E_{k+i+1}[(m+1)r+l+j-k-i-1]},$

2)for $k+i< r+j$ and $i\leq l+j\leq k+i-1$, we have
$X_{E_i[k]}X_{E_j[mr+l]}=X_{E_j[mr+k+i-j]}X_{E_i[l+j-i]}+X_{E_j[mr+i-j-1]}X_{E_{l+j+1}[k+i-l-j-1]},$

3)for other conditions, we have
$X_{E_i[k]}X_{E_j[mr+l]}=X_{E_i[k]\oplus E_j[mr+l]}$.

\nd (2)When $j> i$, then

1)for $k\geq j-i,$ we have
$X_{E_i[k]}X_{E_j[mr+l]}=X_{E_i[j-i-1]}X_{E_{k+i+1}[mr+l+j-k-i-1]}+X_{E_i[mr+l+j-i]}X_{E_j[k+i-j]},$

2)for $k< j-i$ and $i\leq l+j-r\leq k+i-1$, we have
$X_{E_i[k]}X_{E_j[mr+l]}=X_{E_j[(m+1)r+k+i-j]}X_{E_i[l+j-r-i]}+X_{E_j[(m+1)r+i-j-1]}X_{E_{l+j+1}[k+r+i-l-j-1]},$

3)for other conditions, we have
$X_{E_i[k]}X_{E_j[mr+l]}=X_{E_i[k]\oplus E_j[mr+l]}.$
\end{Thm}

\begin{proof}
We only prove (1) and (2) is totally similar to (1).

1) When $k=1,$ by $k+i\geq r+j$ and $1\leq j\leq i\leq
r\Longrightarrow i=r$ and $j=1.$\\
Then by the cluster multiplication theorem in Theorem \ref{XX} or
Theorem \ref{CK}, we have
$$X_{E_r}X_{E_1[mr+l]}=X_{E_r[mr+l+1]}+X_{E_2[mr+l-1]}.$$
 When $k=2,$ by $k+i\geq r+j$ and $1\leq j\leq i\leq
r\Longrightarrow
i=r\ or\ i=r-1.$\\
For $i=r\Longrightarrow j=1\ or\ j=2$:\\
The case for $i=r$ and $j=1$, we have
\begin{eqnarray}
  X_{E_r[2]}X_{E_1[mr+l]} &=& (X_{E_r}X_{E_1}-1)X_{E_1[mr+l]} \nonumber\\
  &=& X_{E_1}(X_{E_r[mr+l+1]}+X_{E_2[mr+l-1]})-X_{E_1[mr+l]} \nonumber\\
  &=& X_{E_1}X_{E_r[mr+l+1]}+(X_{E_1[mr+l]}+X_{E_3[mr+l-2]})-X_{E_1[mr+l]} \nonumber\\
  &=& X_{E_1}X_{E_r[mr+l+1]}+X_{E_3[mr+l-2]}.\nonumber
\end{eqnarray}
The case for $i=r$ and $j=2$, we have
\begin{eqnarray}
  X_{E_r[2]}X_{E_2[mr+l]} &=& (X_{E_r}X_{E_1}-1)X_{E_2[mr+l]} \nonumber\\
  &=& X_{E_r}(X_{E_1[mr+l+1]}+X_{E_3[mr+l-1]})-X_{E_2[mr+l]} \nonumber\\
  &=& X_{E_r[mr+l+2]}+(X_{E_2[mr+l]}+X_{E_r}X_{E_3[mr+l-1]})-X_{E_2[mr+l]} \nonumber\\
  &=& X_{E_r[mr+l+2]}+X_{E_r}X_{E_3[mr+l-1]}.\nonumber
\end{eqnarray}

For $i=r-1\Longrightarrow j=1$:
\begin{eqnarray}
  X_{E_{r-1}[2]}X_{E_1[mr+l]}&=& (X_{E_{r-1}}X_{E_r}-1)X_{E_1[mr+l]} \nonumber\\
  &=& X_{E_{r-1}}(X_{E_r[mr+l+1]}+X_{E_2[mr+l-1]})-X_{E_1[mr+l]} \nonumber\\
  &=& (X_{E_{r-1}[mr+l+2]}+X_{E_1[mr+l]})+X_{E_{r-1}}X_{E_2[mr+l-1]}-X_{E_1[mr+l]} \nonumber\\
  &=& X_{E_{r-1}[mr+l+2]}+X_{E_{r-1}}X_{E_2[mr+l-1]}.\nonumber
\end{eqnarray}\\
Now, suppose it holds for $k\leq n,$ then by induction we have
$$\hspace{-8.0cm}X_{E_i[n+1]}X_{E_j[mr+l]}$$
$$\hspace{-6.8cm}=(X_{E_i[n]}X_{E_{i+n}}-X_{E_i[n-1]})X_{E_j[mr+l]}$$
$$\hspace{-5.3cm}=X_{E_{i+n}}(X_{E_i[n]}X_{E_j[mr+l]})-X_{E_i[n-1]}X_{E_j[mr+l]}$$
$$=X_{E_{i+n}}(X_{E_i[(m+1)r+l+j-i]}X_{E_j[n+i-r-j]}+X_{E_i[r+j-i-1]}X_{E_{n+i+1}[(m+1)r+l+j-n-i-1]})$$
$$\hspace{-0.2cm}-(X_{E_i[(m+1)r+l+j-i]}X_{E_j[n+i-r-j-1]}+X_{E_i[r+j-i-1]}X_{E_{n+i}[(m+1)r+l+j-n-i]})$$
$$\hspace{-4.1cm}=X_{E_i[(m+1)r+l+j-i]}(X_{E_j[n+i+1-r-j]}+X_{E_j[n+i-r-j-1]})$$
$$\hspace{-1cm}+X_{E_i[r+j-i-1]}(X_{E_{n+i}[(m+1)r+l+j-n-i]}
+X_{E_{n+i+2}[(m+1)r+l+j-n-i-2]})$$
$$-(X_{E_i[(m+1)r+l+j-i]}X_{E_j[n+i-r-j-1]}+X_{E_i[r+j-i-1]}X_{E_{n+i}[(m+1)r+l+j-n-i]})$$
$$\hspace{-0.5cm}=X_{E_i[(m+1)r+l+j-i]}X_{E_j[n+i+1-r-j]}+X_{E_i[r+j-i-1]}X_{E_{n+i+2}[(m+1)r+l+j-n-i-2]}.$$

2)  When $k=1,$ by $i\leq l+j\leq k+i-1\Longrightarrow i\leq
l+j\leq i\Longrightarrow i=l+j.$\\
Then by the cluster multiplication theorem in Theorem \ref{XX} or
Theorem \ref{CK}, we have
$$X_{E_{i}}X_{E_j[mr+l]}=X_{E_{l+j}}X_{E_j[mr+l]}=X_{E_j[mr+l+1]}+X_{E_j[mr+l-1]}$$
  When $k=2,$ by $i\leq l+j\leq k+i-1\Longrightarrow i\leq
l+j\leq i+1\Longrightarrow i=l+j\ or\ i+1=l+j$:\\
For $i=l+j$, we have
$$\hspace{-2.0cm}X_{E_{i}[2]}X_{E_j[mr+l]}=X_{E_{l+j}[2]}X_{E_j[mr+l]}=(X_{E_{l+j}}X_{E_{l+j+1}}-1)X_{E_j[mr+l]}$$
$$\hspace{0.5cm}=(X_{E_j[mr+l+1]}+X_{E_j[mr+l-1]})X_{E_{l+j+1}}-X_{E_j[mr+l]}$$
$$\hspace{2.1cm}=X_{E_j[mr+l+2]}+X_{E_j[mr+l]}+X_{E_{l+j+1}}X_{E_j[mr+l-1]}-X_{E_j[mr+l]}$$
$$\hspace{-1.6cm}=X_{E_j[mr+j+2]}+X_{E_{l+j+1}}X_{E_j[mr+l-1]}.$$
For $i+1=l+j$, we have
$$\hspace{-1.6cm}X_{E_{i}[2]}X_{E_j[mr+l]}=X_{E_{l+j-1}[2]}X_{E_j[mr+l]}=(X_{E_{l+j-1}}X_{E_{l+j}}-1)X_{E_j[mr+l]}$$
$$\hspace{0.5cm}=(X_{E_j[mr+l+1]}+X_{E_j[mr+l-1]})X_{E_{l+j-1}}-X_{E_j[mr+l]}$$
$$\hspace{2.5cm}=X_{E_j[mr+l+1]}X_{E_{l+j-1}}+(X_{E_j[mr+l]}+X_{E_j[mr+l-2]})-X_{E_j[mr+l]}$$
$$\hspace{-1.5cm}=X_{E_j[mr+l+1]}X_{E_{l+j-1}}+X_{E_j[mr+l-2]}.$$\\
Suppose it holds for $k\leq n,$ then by induction we have
$$\hspace{-8.0cm}X_{E_i[n+1]}X_{E_j[mr+l]}$$
$$\hspace{-6.5cm}=(X_{E_i[n]}X_{E_{i+n}}-X_{E_i[n-1]})X_{E_j[mr+l]}$$
$$\hspace{-5.0cm}=(X_{E_i[n]}X_{E_j[mr+l]})X_{E_{i+n}}-X_{E_i[n-1]}X_{E_j[mr+l]}$$
$$\hspace{-1.3cm}=(X_{E_j[mr+n+i-j]}X_{E_i[l+j-i]}+X_{E_j[mr+i-j-1]}X_{E_{l+j+1}[n+i-l-j-1]})X_{E_{i+n}}$$
$$\hspace{-1.4cm}-(X_{E_j[mr+n+i-j-1]}X_{E_i[l+j-i]}+X_{E_j[mr+i-j-1]}X_{E_{l+j+1}[n+i-l-j-2]})$$
$$\hspace{-4.45cm}=(X_{E_j[mr+n+i+1-j]}+X_{E_j[mr+n+i-j-1]})X_{E_i[l+j-i]}$$
$$\hspace{-2.9cm}+(X_{E_{l+j+1}[n+i-l-j]}+X_{E_{l+j+1}[n+i-l-j-2]})X_{E_j[mr+i-j-1]}$$
$$\hspace{-1.3cm}-(X_{E_j[mr+n+i-j-1]}X_{E_i[l+j-i]}+X_{E_j[mr+i-j-1]}X_{E_{l+j+1}[n+i-l-j-2]})$$
$$\hspace{-2.4cm}=X_{E_j[mr+n+i+1-j]}X_{E_i[l+j-i]}+X_{E_j[mr+i-j-1]}X_{E_{l+j+1}[n+i-l-j]}.$$

3)  It is trivial by the definition of the  Caldero-Chapton map.
\end{proof}

Note that in the above section, we have already obtained a
$\mathbb{Z}$-basis for the cluster algebra of
$\widetilde{A}_{n,n},\widetilde{D}\ or\ \widetilde{E}$. Now by
using Theorem \ref{16} and Proposition \ref{17} we can easily
express $X_{E_i[m]}X_{E_j[n]}$ as a $\mathbb{Z}-$combinatorics of
the basis for the cluster algebra of
$\widetilde{A}_{n,n},\widetilde{D}\ or\ \widetilde{E}$ where
$E_i[m]\ and\ E_j[n]$ are regular modules in one fixed tube. In
the following section, we will explain it by an example.

\section{an example}
We consider a tube with rank 3. By Theorem \ref{16}, we can easily
obtain the following proposition.
\begin{Prop}\label{18}
(1) For $n\geq 3m+1$, then
$$X_{E_2[3m+1]}X_{E_1[n]}=X_{E_2}X_{E_1[n+3m]}+X_{E_1[n+3m-3]}+X_{E_2}X_{E_1[n+3m-6]}+X_{E_1[n+3m-9]}$$
$$\hspace{2.3cm}+\cdots+X_{E_2}X_{E_1[n-3m+6]}+X_{E_1[n-3m+3]}
+X_{E_2}X_{E_1[n-3m]},$$

(2) For $n\geq 3m+2$, then
$$X_{E_2[3m+2]}X_{E_1[n]}=X_{E_2[n+3m+2]}+X_{E_2}X_{E_2[n+3m-1]}+X_{E_2[n+3m-4]}+X_{E_2}X_{E_2[n+3m-7]}$$
$$\hspace{-0.3cm}+\cdots+X_{E_2[n-3m+2]}+X_{E_2}X_{E_2[n-3m-1]},$$

(3) For $n\geq 3m+3$, then
$$\hspace{-0.56cm}X_{E_2[3m+3]}X_{E_1[n]}=X_{E_1[n+3m+3]}+X_{E_3[n+3m+1]}+X_{E_2[n+3m-1]}+X_{E_1[n+3m-3]}$$
$$\hspace{1.1cm}+X_{E_3[n+3m-5]}+X_{E_2[n+3m-7]}+\cdots+X_{E_1[n-3m+3]}$$
$$\hspace{0.35cm}+X_{E_3[n-3m+1]}+X_{E_2[n-3m-1]}+X_{E_1[n-3m-3]}.$$
\end{Prop}
\begin{proof}
(1)  By Theorem \ref{16}, then
$$\hspace{-8cm}X_{E_2[3m+1]}X_{E_1[n]}$$
$$\hspace{-6.6cm}=X_{E_2[n+2]}X_{E_1[3m-1]}+X_{E_2}X_{E_1[n-3m]}$$
$$\hspace{-4.2cm}=X_{E_2[3m-2]}X_{E_1[n+3]}+X_{E_1[n-3m+3]}+X_{E_2}X_{E_1[n-3m]}$$
$$\hspace{-1.2cm}=X_{E_2[n+5]}X_{E_1[3m-4]}+X_{E_2}X_{E_1[n-3m+6]}+X_{E_1[n-3m+3]}
+X_{E_2}X_{E_1[n-3m]}$$
$$\hspace{0.05cm}=X_{E_2[3m-5]}X_{E_1[n+6]}+X_{E_1[n-3m+9]}+X_{E_2}X_{E_1[n-3m+6]}+X_{E_1[n-3m+3]}
+X_{E_2}X_{E_1[n-3m]}$$
$$\vdots$$
$$\hspace{-2.0cm}=X_{E_2}X_{E_1[n+3m]}+X_{E_1[n+3m-3]}+X_{E_2}X_{E_1[n+3m-6]}+X_{E_1[n+3m-9]}$$
$$\hspace{-2.6cm}+\cdots+X_{E_2}X_{E_1[n-3m+6]}+X_{E_1[n-3m+3]}
+X_{E_2}X_{E_1[n-3m]}.$$

(2) By  Theorem \ref{16}, then
$$\hspace{-8cm}X_{E_2[3m+2]}X_{E_1[n]}$$
$$\hspace{-5.7cm}=X_{E_2[n+2]}X_{E_1[3m]}+X_{E_2}X_{E_2[n-3m-1]}$$
$$\hspace{-2.9cm}=X_{E_2[3m-1]}X_{E_1[n+3]}+X_{E_2[n-3m+2]}+X_{E_2}X_{E_2[n-3m-1]}$$
$$\hspace{0.1cm}=X_{E_2[n+5]}X_{E_1[3m-3]}+X_{E_2}X_{E_2[n-3m+5]}+X_{E_2[n-3m+2]}+X_{E_2}X_{E_2[n-3m-1]}$$
$$\hspace{0.5cm}=X_{E_2[3m-4]}X_{E_1[n+6]}+X_{E_2[n-3m+8]}+X_{E_2}X_{E_2[n-3m+5]}+X_{E_2[n-3m+2]}+X_{E_2}X_{E_2[n-3m-1]}$$
$$\vdots$$
$$\hspace{-0.8cm}=X_{E_2[n+3m+2]}+X_{E_2}X_{E_2[n+3m-1]}+X_{E_2[n+3m-4]}+X_{E_2}X_{E_2[n+3m-7]}$$
$$\hspace{-4.5cm}+\cdots+X_{E_2[n-3m+2]}+X_{E_2}X_{E_2[n-3m-1]}.$$

(3)  By  Theorem \ref{16}, then
$$\hspace{-8cm}X_{E_2[3m+3]}X_{E_1[n]}$$
$$\hspace{-5.0cm}=X_{E_2[n+2]}X_{E_1[3m+1]}+X_{E_2}X_{E_3[n-3m-2]}$$
$$\hspace{-3.3cm}=X_{E_2[n+2]}X_{E_1[3m+1]}+X_{E_2[n-3m-1]}+X_{E_1[n-3m-3]}$$
$$\hspace{-1.26cm}=X_{E_2[3m]}X_{E_1[n+3]}+X_{E_3[n-3m+1]}+X_{E_2[n-3m-1]}+X_{E_1[n-3m-3]}$$
$$\vdots$$
$$\hspace{-4.3cm}=X_{E_1[n+3m+3]}+X_{E_3[n+3m+1]}+X_{E_2[n+3m-1]}$$
$$\hspace{-3.4cm}+X_{E_1[n+3m-3]}
+X_{E_3[n+3m-5]}+X_{E_2[n+3m-7]}$$
$$+\cdots+X_{E_1[n-3m+3]}+X_{E_3[n-3m+1]}+X_{E_2[n-3m-1]}+X_{E_1[n-3m-3]}.$$
\end{proof}
\begin{Cor}\label{20}
When $n=3k+1$,we can rewrite Proposition \ref{18} as following:

(1) For $n\geq 3m+1$, then
$$\hspace{-1.8cm}X_{E_2[3m+1]}X_{E_1[n]}=X_{E_1[n+3m+1]}+X_{E_1[n+3m-1]}+X_{E_1[n+3m-3]}+\cdots$$
$$\hspace{2.4cm}+X_{E_1[n-3m+5]}+X_{E_1[n-3m+3]}+X_{E_1[n-3m+1]}
+X_{E_1[n-3m-1]},$$

(2) For $n\geq 3m+2$, then
$$\hspace{-0.6cm}X_{E_2[3m+2]}X_{E_1[n]}=X_{E_2[n+3m+2]}+X_{E_2[n+3m]}+X_{E_2[n+3m-2]}+X_{E_2[n+3m-4]}$$
$$\hspace{0.6cm}+\cdots+X_{E_2[n-3m+2]}+X_{E_2[n-3m]}+X_{E_2[n-3m-2]},$$

(3) For $n\geq 3m+3$, then
$$\hspace{-2.7cm}X_{E_2[3m+3]}X_{E_1[n]}=X_{E_1[n+3m+3]}+X_{E_3[n+3m+1]}+X_{E_2[n+3m-1]}$$
$$\hspace{0.8cm}+X_{E_1[n+3m-3]}
+X_{E_3[n+3m-5]}+X_{E_2[n+3m-7]}+\cdots$$
$$\hspace{2.4cm}+X_{E_1[n-3m+3]}+X_{E_3[n-3m+1]}+X_{E_2[n-3m-1]}+X_{E_1[n-3m-3]}.$$
\end{Cor}
\begin{proof}
When $n=3k+1$,we have the following equations:
$$X_{E_2}X_{E_1[n]}=X_{E_1[n+1]}+X_{E_1[n-1]}$$
$$X_{E_2}X_{E_2[n-1]}=X_{E_2[n]}+X_{E_2[n-2]}$$
Then the proof is immediately finished from Proposition \ref{18}.
\end{proof}

In the same way by using Theorem \ref{16}, we have the following
proposition.
\begin{Prop}\label{19}
(1) For $n\geq 3m+1$, then
$$X_{E_1[3m+1]}X_{E_1[n]}=X_{E_1}X_{E_1[n+3m]}+X_{E_1[n+3m-3]}+X_{E_1}X_{E_1[n+3m-6]}+X_{E_1[n+3m-9]}$$
$$\hspace{1.7cm}+\cdots+X_{E_1}X_{E_1[n-3m+6]}+X_{E_1[n-3m+3]}
+X_{E_1}X_{E_1[n-3m]},$$

(2) For $n\geq 3m+2$, then
$$\hspace{-1.5cm}X_{E_1[3m+2]}X_{E_1[n]}=X_{E_1[2]}(X_{E_1[n+3m]}+X_{E_1[n+3m-6]}+\cdots+X_{E_1[n-3m]}),$$

(3) For $n\geq 3m+3$, then
$$X_{E_1[3m+3]}X_{E_1[n]}=X_{E_1[n+3m+3]}+X_{E_2}X_{E_1[n+3m]}+X_{E_1[n+3m-3]}+X_{E_2}X_{E_1[n+3m-6]}$$
$$\hspace{1.0cm}+\cdots+X_{E_1[n-3m+3]}+X_{E_2}X_{E_1[n-3m]}+X_{E_1[n-3m-3]}.$$
\end{Prop}
\begin{Remark}
By the same method in Corollary \ref{20}, we can also rewrite
Proposition \ref{19} if we consider these different n. Here we
omit it.
\end{Remark}
\begin{Remark}
In fact, we can also check Proposition \ref{18} Corollary \ref{20}
and Proposition \ref{19} by induction.
\end{Remark}
If we consider the quiver is one of
$\widetilde{A}_{n,n},\widetilde{D}\ and\ \widetilde{E}$. Then from
the above discussions and Proposition \ref{17}, we can easily
express $X_{E_i[m]}X_{E_j[n]}$ as a $\mathbb{Z}-$combinatorics of
the basis in this tube with rank 3.

\section{A $\mathbb{Z}$-basis for the cluster algebra of affine type}
In \cite{Dupont1} Theorem 4.21, G. Dupont asserts that the
following set is a $\mathbb{Z}$-basis for the cluster algebra of
type $\widetilde{A}$:
$$\{cluster\  monomials\}\cup \{X_{M_{\lambda}^{\oplus n}\oplus E}:n\geq 1,E\ are\ rigid\ regular\ modules\}.$$
We will use Theorem \ref{7} and the above result in \cite{Dupont1}
to give a $\mathbb{Z}$-basis for the cluster algebra of affine
type.

\nd \emph{Proof of Theorem \ref{30}.} By Theorem \ref{60}, we only
prove that the set $\mathcal{B}(Q')$ is a $\mathbb{Z}$-basis for
the cluster algebra of $\widetilde{A}_{r,s}$ for $r\neq s$. By
Theorem 3.25 in \cite{Dupont1} and Theorem \ref{16}, we have
\begin{eqnarray}
  X_{M_{\lambda}^{\oplus n}\oplus E} &=& (X_{\delta_{i,1}}-X_{q.rad M_{E_{i,1}}/E_{i,1}})^{n}X_{E} \nonumber\\
  &=& (X_{\delta_{i,1}}^{n}+\sum_{\mathrm{\underline{dim}}(T'\oplus
R')\prec n\delta} a_{T'\oplus R'}X_{T'\oplus
R'})X_{E} \nonumber\\
  &=& X_{T\oplus R}+\sum_{\mathrm{\underline{dim}}(T''\oplus
R'')\prec \mathrm{\underline{dim}}(T\oplus R)} a_{T''\oplus
R''}X_{T''\oplus R''}\nonumber
\end{eqnarray}
where $a_{T'\oplus R'},a_{T''\oplus R''}\in \mathbb{Z}$,
$\mathrm{\underline{dim}}(T\oplus
R)=\mathrm{\underline{dim}}(M_{\lambda}^{\oplus n}\oplus E)$ and
$X_{T\oplus R},X_{T''\oplus R''}\in \mathcal{B}(Q')$. Thus from
the above upper triangle equations and these dimension vectors of
$M_{\lambda}^{\oplus n}\oplus E$ in \cite{Dupont1} are different,
we have
$$\{X_{M_{\lambda}^{\oplus n}\oplus E}:n\geq 1,E\ are\ rigid\
regular\ modules\}=\{X_{T\oplus
R}|\mathrm{\underline{dim}}(T_1\oplus
R_1)\neq\mathrm{\underline{dim}}(T_2\oplus R_2)\}$$ where $R$ is 0
or any regular exceptional module, $T$ is any indecomposable
regular module with self-extension and there are no extension
between $R$ and $T$.

Also it is easy to see that $$\{cluster\
monomials\}=\{X_{L}|Ext^{1}_{\mc(Q)}(L,L)=0\}$$ Therefore
$\mathcal{B}(Q')$ is a $\mathbb{Z}$-basis for the cluster algebra
of $\widetilde{A}_{r,s}$.

\end{document}